\documentclass[10pt]{amsart}

\usepackage{amsfonts}
\usepackage{amssymb}
\usepackage{amsmath, amsthm, latexsym}
\usepackage[all]{xy}
\usepackage{hyperref}
\usepackage[margin=1.5in]{geometry}
\usepackage{color}

\def \c{\mathbb{C}}

\def \T{\mathbb{T}}
\def\B{\mathcal B}
\def\N{\mathcal N}
\def\W{\mathcal W}

\def \t{\mathfrak{t}}
\def \g{\mathfrak{g}}
\def\s{\mathfrak{s}}
\def\m{\mathfrak{l}}
\def\R{\c[\s]}

\def\coh{cohomology}
\def\eq{equivariant}
\def\ec{equivariant cohomology}
\def\crm{cohomology restriction map}

\def\ef{equivariantly formal}

\def\T{\times}

\theoremstyle{plain}
\newtheorem{Th}{Theorem}[section]
\newtheorem{Lem}[Th]{Lemma}

\newtheorem{Cor}[Th]{Corollary}

\theoremstyle{definition}

\newtheorem{Rem}[Th]{Remark}

\begin{document}

\title{Springer's Weyl Group Representation  via Localization}

\author{Jim Carrell}
\address{Department of Mathematics, University of British Columbia, Vancouver, B.C., Canada}
\email{carrell@math.ubc.ca}

\author{Kiumars Kaveh}
\address{Department of Mathematics, University of Pittsburgh,
Pittsburgh, PA, USA.}
\email{kaveh@pitt.edu}

\noindent
\date{\today}

\maketitle

\begin{abstract} Let $G$ denote a reductive algebraic group over $\c$
and $x$ a nilpotent element of its Lie algebra $\g$. The Springer variety $\B_x$
is the closed subvariety of the flag variety $\B$ of $G$ parameterizing the 
Borel subalgebras of $\g$ containing $x$. It has the remarkable property that 
the Weyl group $W$ of $G$ admits a representation on the \coh\ of $\B_x$
even though $W$ rarely acts on $\B_x$ itself. Well-known constructions of this action 
due to Springer et al use technical machinery from algebraic geometry. 
The purpose of this note is to describe an elementary approach that gives this action
when $x$ is what we call parabolic-surjective. The idea is to 
use localization to construct an action of $W$ on
the \ec\ algebra $H_S^*(\B_x)$, where $S$ is a certain algebraic subtorus of $G$.
This action descends to $H^*(\B_x)$ via the forgetful map and
gives the desired representation. The parabolic-surjective case includes all nilpotents of type $A$ and,
more generally, all nilpotents for which it is known that $W$ acts on 
$H_S^*(\B_x)$ for some torus $S$ (see \cite{KP} and \cite{GM}). 
Our result is deduced from a general theorem describing  
when a group action on the \coh\ of the fixed point set of a torus action
on a space lifts to the full \coh\ algebra of the space. 
\end{abstract}

\section{Introduction}
Let $G$ be a reductive linear algebraic group over $\c$
with Lie algebra $\g$, and fix a  
maximal torus $T$ and a  Borel subgroup $B$ of $G$ such that $T\subset B$. 
The flag variety $\B=G/B$ of $G$ will be viewed as the variety of all Borel subgroups of $G$
or, equivalently, as the variety of all Borel subalgebras of  $\g$.
Let $W=N_G(T)/T$ be the Weyl group  of the pair $(G,T)$, and
recall that $W$ acts topologically on $\B$. Thus the \coh\
algebra $H^*(\B)$ admits a representation as a graded $W$-module,
which is well-known to be isomorphic with the graded $W$-algebra
$ \c[\t]/I_W^+$, the coinvariant algebra of $W$. Here $\t$ is the Lie algebra of $T$ 
and $I_W^+$ is the ideal in $\c[\t]$ generated by the nonconstant homogeneous $W$-invariants. 
Note: throughout this paper, $H^*(Y)$ will denote the standard cohomology algebra
of a space $Y$ with complex coefficients.

A celebrated theorem of T. A. Springer \cite{SPRINGER1,SPRINGER2} 
says that if $x$ is a nilpotent element of 
$\g$ and $\B_x$ is the  Springer variety associated to $x$,
namely the closed subvariety of $\B$  consisting of all
Borel subalgebras of $\g$ containing $x$,
 then there is a graded $\c$-algebra representation of $W$  on  $H^*(\B_x)$ 
so that the \crm\ $i_x^*: H^*(\B)\to H^*(\B_x)$ associated to
the inclusion $i_x:\B_x \hookrightarrow \B$ is $W$-\eq\ (see \cite{HS} for the proof of 
$W$-equivariance). As remarked often, the existence of this  
representation is quite surprising because $W$ itself does
not usually act on $\B_x$, exceptions being
when $x=0$ or $x$ is regular in $\g$, and the usual definitions
require  a lot of technical machinery. Subsequent 
definitions involve either replacing 
$H^*(\B_x)$ with an isomorphic algebra on which $W$ is known to act
(cf. \cite{KRAFT,DP,COMP}) or replacing $\B_x$ with a space
having both a $W$-action and isomorphic \coh\ algebra 
(cf. \cite{Ross,Slo,Treu}).

Our plan is to give a simple new 
construction of Springer's representation
when $x$ is a nilpotent element of $\g$ 
which is what we call {\em parabolic-surjective}.
We call a nilpotent $x\in \g$  {\em parabolic} if it is regular 
in a Levi subalgebra $\m$ of $\g$, and we say
$x$ is {\em parabolic-surjective} if in addition
the \crm\ $i_x^*$ is surjective. 
The key idea is to exploit the fact that in the parabolic-surjective case, 
there exists subtorus $S$ of $G$ acting on $\B_x$ 
so that $(\B_x)^S$ is stable under $W$ 
so localization and the parabolic-surjective condition can be used to obtain
an $\R$-module  action of $W$ on the torus \ec\ $H^*_S(\B_x)$ descending to 
an  action on $H^*(\B_x)$ which turns out to coincide with Springer's action. 

Recently, the existence of this $W$-action on $H^*_S(\B_x)$
was established by finding a geometric realization of Spec$(H_S^*(\B_n))$:
see \cite{KP} for the general parabolic-surjective setting and \cite{GM}
for $x$ of type $A$. We note that although  \cite{GM} only treats the 
type $A$ case, their argument is valid for all parabolic-surjective $x\in \g$
fo all $\g$. The paper \cite{AH} establishes this action by 
employing Tanisaki's presentation of $H^*(\B_x)$ in type $A$. 

We now state the main result which will give Springer's action
for parabolic-surjective $x$ by first obtaining it for $H^*_S(\B_x)$.
\begin{Th} \label{MAINRESULT} 
Let $Y$ be a projective variety with vanishing odd cohomology. Also suppose that we have actions of an algebraic 
torus $S$ and a finite group $\W$ on $Y$ and these two actions commute with each other.
Let $X$ be an $S$-stable subvariety of $Y$ 
such that the \crm\ $i^*:H^*(Y)\to H^*(X)$ is surjective, where $i:X \hookrightarrow Y$ is the inclusion. 
Then if  $\W$ acts on $H^*(X^S)$ so that the \crm\ $(i_S)^*:H^*(Y^S)\to H^*(X^S)$
induced by the inclusion $i_S:X^S \hookrightarrow Y^S$ 
is $\W$-\eq, then $\W$ also acts on $H_S^*(X)$ by 
graded $\W$-algebra $\c[\s]$-module isomorphisms. 
Moreover, the natural map $H_S^*(X)\to H^*(X)$
induces a representation of $\W$ on the graded algebra $H^*(X)$ 
compatible with its action on $H^*(Y)$ such that all the maps in the 
following commutative diagram  are $\W$-\eq:
\begin{equation} \label{equ-comm-diag3}
\xymatrix{
H^*_S(Y) \ar[r]^{\iota^*} \ar[d] & H^*_S(X) \ar[d] \\
H^*(Y) \ar[r]^{i^*} & H^*(X). \\
}
\end{equation}

\end{Th}

In the above diagram,  ${\iota^*}$ is the map on \ec\ induced by $i:X\hookrightarrow Y$. The proof, 
given in the next section, is an application of the localization theorem. 
The reader can easily reformulate this result as a statement involving 
topological torus actions. 
In the final section, we will verify the above assertions about  parabolic-surjective  
Springer varieties.

\medskip
{\bf Acknowledgement} We would like to thank Shrawan Kumar for  useful comments and Megumi Harada for pointing out
the paper \cite{AH}.

\section{proof of the main theorem}

We will begin by reviewing some facts about \ec. Excellent references for 
the facts below are \cite{AB,BRION}.
Recall that all \coh\ is over $\c$. 
Let $Y$ be a complex projective variety with vanishing odd \coh\
admitting a nontrivial action $(S,Y)$ by an algebraic torus $S\cong (\c^*)^\ell$. In particular,
the fixed point set $Y^S$ is nontrivial. Recall that the $S$-\eq\ \coh\ algebra $H^*_S(Y)$ of $Y$  
is defined as the \coh\ algebra
$H^*(Y_S)$  of the Borel space $Y_S=(Y\T E)/S$,
where $E$ is a contractible space with a free $S$-action and  $S$ acts diagonally
on the product. The projection $Y \times E \to E$ induces a map $Y_S = (Y \T E) / S \to E/S$ 
which in turn induces an $H^*(E/S)$-module structure on $H_S^*(Y)$.
On the other hand, the inclusion $\nu_Y: Y \hookrightarrow Y_S$ along a fibre gives a 
map from $H^*_S(Y)$ to the ordinary cohomology $H^*(Y)$. 
Moreover, there is a natural identification $H^*(E/S) \cong \c[\s]$, where $\s=\mathrm{Lie}(S)$.
Note that by the Kunneth formula, $H^*_S(Y^S)=\c[\s]\otimes H^*(Y^S)$.
When $H^*_S(Y)$ is a free $\c[\s]$-module,  the action $(S,Y)$ is said to be {\em \ef}. 
It is well-known that equivariant formality is implied by the vanishing of odd cohomology
of $Y$. The proof of our main result is based on the following well-known result, the
first assertion of which is a special case of the localization theorem.
\begin{Th}\label{LOCTH}
If $(S,Y)$ is \ef, then the inclusion mapping
$j_Y:Y^S \hookrightarrow Y$ induces an injection $(j_Y)_S^*: H^*_S(Y)\to H^*_S(Y^S)$. 
Moreover, the map $(\nu_Y)^*$ fits into an exact sequence
\begin{equation}\label{EXSEX}
0\to \R^+H^*_S(Y)\to H^*_S(Y)\xrightarrow{(\nu_Y)^*} H^*(Y)\to 0,
\end{equation}
where $\R^+$ is the augmentation ideal, generated by all the nonconstant homogeneous polynomials. 
\end{Th}

We now prove Theorem \ref{MAINRESULT}. 
Assume  that $\W$ is a finite group acting on $Y$ such that the action $(\W,Y)$
commutes with $(S,Y)$. Then $\W$ acts linearly on both $H^*(Y)$ and $H^*(Y^S)$. 
Furthermore, it acts on $H^*_S(Y)$ and $H^*_S(Y^S)$ as 
$\R$-module isomorphisms so that the  map $j_Y^*:H^*_S(Y)\to H^*_S(Y^S)$
induced by the inclusion $j_Y:Y^S \hookrightarrow Y$ is a $\W$-\eq\ $\R$-module injection. 
Let $X$ be an $S$-stable subvariety of $Y$ such that the \crm\ $i^*:H^*(Y)\to H^*(X)$ is surjective.

Since $\W$ acts on $H^*(X^S)$ and the restriction map $H^*(Y^S)\to H^*(X^S)$ is $\W$-\eq, 
$\W$ also acts on $H^*_S(X^S)$
as a group of  $\R$-module isomorphisms so that the natural map $\mu:H^*_S(Y^S)\to H^*_S(X^S)$ 
is a $\W$-\eq\ $\R$-module homomorphism. 
Now consider the commutative diagram
 \begin{equation} \label{equ-comm-diag2}
\xymatrix{
H^*_S(Y) \ar[r]^{\iota^*} \ar[d]^{j^*_{Y}} & H^*_S(X) \ar[d]^{j^*_X} \\
H^*_S(Y^S) \ar[r]^{\mu} & H_S^*(X^S) \\
}
\end{equation}
We will define the action of $\W$ on $H^*_S(X)$ by imposing the requirement that
$\iota^*$ be a $\W$-module homomorphism. To show this action  is well-defined,
it suffices to show that the kernel of $\iota^*$ is a $\W$-submodule. Suppose then that 
$\iota^*(a)=0$. By assumption, $H^*(Y) \to H^*(X)$ is surjective and thus $X$ also has vanishing odd cohomology. It follows that $(S,X)$ is equivariantly 
formal and hence $j_X^*$ is injective.
Thus to show $\iota^*(w\cdot a)=0$ for any $w\in \W$, it suffices to show that 
 $j^*_X \iota^*(w\cdot a)=0$. But
$$j^*_X \iota^*(w\cdot a)=\mu j_{Y}^*(w\cdot a)=w\cdot  \mu j_{Y}^*(a)=wj^*_X \iota ^*(a)=0,$$
since $\mu j_{Y}^*$ is a $\W$-module homomorphism.
Thus, $\W$ acts on $H_S^*(X)$ as claimed. 
It  follows from this argument that $j_X^*$ is $\W$-\eq.
To show that $\W$ acts on $H^*(X)$,
consider the exact sequence (\ref{EXSEX}) for $X.$ 
As above, we  may define the $\W$-action by requiring that $(\nu_X)^*$ be \eq.
It suffices to show its kernel is $\W$-stable.
But if $(\nu_X)^*(a)=0$, then $a=fb$ for some $f\in \c[\s]^+$
and $b\in H^*_S(X)$. Thus,
$$w\cdot a=w\cdot fb=f(w\cdot b) \in \c[\s]^+H^*_S(X)=\ker(\nu_X)^*.$$
Finally, we remark that the above definitions make the diagram 
(\ref{equ-comm-diag3}) commutative.
\qed

 \medskip
If one omits the assumption that  \crm\ $i^*:H^* (Y)\to H^*(X)$ is surjective,
the best one can hope for is that $\W$ acts on the image $i^*(H^*(Y))$.
The following result gives a sufficient condition for $\W$ to act in
this case.
\begin{Th} \label{COR} \label{COR1} Assume that the setup in Theorem  \ref{MAINRESULT} 
holds except for the assumption that $i^*:H^* (Y)\to H^*(X)$ is surjective,
and, in addition, assume also that $H^*(X)$ has vanishing odd \coh. 
Then there exists an action of $\W$ on $\iota^*(H_S^*(Y))$ by
 $\c[\s]$-module isomorphisms. Moreover, if $\iota^*(H_S^*(Y))$ is free
of rank  $\dim i^*(H^*(Y))$, then the action of $\W$ on $\iota^*(H_S^*(Y))$ descends to $i^*(H^*(Y))$ 
so that the \crm\ $i^*: H^*(Y) \to i^*(H^*(Y))$ is $\W$-\eq.
\end{Th}
\proof For the first assertion, we have to show that the kernel of $\iota^*$
is $\W$-invariant. Since $H^*(X)$ has vanishing odd \coh, 
$j_X^*$ is injective, so this follows from the argument above.
Next, note that if $\N$ denotes a free $\R$-module of finite rank, then
the  $\c$-vector space dimension of $\N / \R^+ \N$ is equal to the rank of $\N$. 
Thus it follows by assumption that the sequence
$$0\to \R^+ \iota^*(H_S^*(Y)) \to \iota^*(H_S^*(Y)) \xrightarrow{\nu} i^*(H^*(Y)) \to 0$$
is exact, where $\nu$ is the restriction of $(\nu_X)^*$. Hence, as above,  the kernel of 
$\iota^*(H_S^*(Y)) \to i^*(H^*(Y))$ is $\W$-stable, so $\W$ acts on $i^*(H^*(Y))$. 
Moreover, the map $H^*(Y) \to i^*(H^*(Y))$ is $\W$-\eq.\qed

\begin{Rem}  In the case when $S = \c^*$, since $(S,X)$ is assumed to be \ef, the module $\iota^*(H_S^*(Y))$ is always free. This is
because $\R$ is a principal ideal domain and $H_S^*(X)$ is free. 
\end{Rem}

\section{The Weyl group action on $H_S^*(\B_x)$}
We now return to the parabolic-surjective setting. First, recall
 that $W$ acts as a group of homeomorphisms of $\B$ which commute with $T$.
Let $K$ be a maximal compact subgroup in $G$ such that $H=K\cap T$ is a maximal torus in $K$. Then
the natural mapping $K/H\to \B$ is a homeomorphism; but
$W=N_K(H)/H$ acts on $K/H$ (from the left) by $w\cdot kH=k\dot{w}^{-1}H$, 
where $\dot{w} \in N_K(H)$ is a representative of $w$.
Thus $W$ acts on $\B$ as asserted. Since this action commutes with the action of $H$ on $K/H$, the group
$W$ acts on both $H^*_T(\B)=H_{H}^*(K/H)$ and $H^*(\B)$
and the natural mapping $H^*_T(\B)\to H^*(\B)$ is $W$-\eq.

In order to apply Theorem \ref{MAINRESULT}, we need the following lemma from \cite{COMP}.

\begin{Lem}\label{LEMMA}  Let $x\in  \g$ be nilpotent, and suppose $x$ is a regular element  in
the Lie algebra of the Levi $L=C_G(S)$ for a subtorus $S$ of $T$. Then $S$ acts on
$\B_x$ with exactly $[W:W_L]$ fixed points.
Moreover, every component of $\B^S$ contains exactly one point of $(\B_x)^S$,
so $W$ acts on $H^*((\B_x)^S)$ so that the \crm\ $H^*(\B^S)\to H^*((\B_x)^S)$ is $W$-\eq\ and surjective.
\end{Lem}

\proof By assumption, $\B_x$ is the variety of Borel subalgebras of $\g$ containing $x$, so 
$\B_x$ is stable under the action of $S$ on $\g$. Each irreducible component of $\B^S$ 
is isomorphic to the flag variety of $L$, so each component contains
a unique fixed point of the  one parameter group exp$(tx)$,
$t\in \c$, since $x$ is regular in $\m$. It follows that the \crm\ $H^*(\B^S)\to H^*((\B_x)^S)$
is surjective. Moreover, since $W$ permutes the components of $\B^S$,
and each component contains a unique point of $\B_x$,
we can uniquely define an action of $W$ on $H^*((\B_x)^S)=H^0((\B_x)^S)$
by requiring that the \crm\ be \eq. Finally, it is well-known that the number of
components of $\B^S$ is the index $[W:W_L]$ of the Weyl group of $L$ in $W$. \qed

\begin{Rem} In fact, $W$ actually acts on $(\B_x)^S$ itself
(see \cite{COMP} or  \cite[Lemma 6.3 and p. 137]{JBC}). 
\end{Rem}

Let us now return to the problem considered in the introduction. As above,
$G$ is reductive linear algebraic group over $\c$ and $\B$ is its flag variety.
By the main result (Theorem \ref{MAINRESULT}), we have
\begin{Cor} Let $x$ be a parabolic-surjective nilpotent in $\g$, say 
$x$ is regular in the Lie algebra of the Levi subgroup $C_G(S)$. 
Then $W$ acts on $H^*_S(\B_x)$ and this action descends to 
to $H^*(\B_x)$ so that the diagram (1) is commutative for $Y=\B$ and $X=\B_x$.
Consequently, this $W$-action is Springer's representation.
\end{Cor}
 
 \proof The only thing to show is that this action of $W$ on $H^*(\B_x)$ coincides with Springer's 
 representation. But this follows  since $i^*:H^*(\B)\to H^*(\B_x)$ 
 is $W$-\eq\ by \cite{HS}. \qed
 
 \medskip
  This seems to give the most elementary construction of Springer's action in type $A$. 
 It was originally conjectured in  \cite{KRAFT} that for any $x\in \s\m(n,\c)$, the action of $W=S_n$ 
 on  $H^*(\B_x)$ is equivalent to the action of $W$ on the coordinate ring
 $A(\t \cap C_y)$ of the schematic intersection of the diagonal matrices in $\s\m(n,\c)$ and
the closure $C_y$  in $\s\m(n,\c)$ of the conjugacy class of
 of the nilpotent $y$ dual to $x$. This was immediately verified in \cite{DP}
 where it was shown that $H^*(\B_x) \cong A(\t \cap C_y)$ as graded $W$-algebras. This 
 isomorphism was extended in \cite{COMP}
 to the case of parabolic-surjective nilpotents 
 in an arbitrary $\g$ which satisfy some additional conditions. 
 Here the  nilpotent $y$ dual to $x$ turns out to be a Richardson element in the 
 nilradical of the parabolic subalgebra of $\g$
 associated to the Levi $\m$ in which $x$ is a regular nilpotent.
 We refer to \cite{COMP} for more details.

Finally, let us mention that by a well-known of result DeConcini, Lusztig
and Procesi \cite{DLP}, $\B_x$ has vanishing odd \coh\ for any nilpotent $x\in \g$. 
Moreover, the  Jacobson-Morosov lemma guarantees that every 
Springer variety  $\B_x$  has a torus action $(S,\B_x)$. This suggests that
$H^*_S(\B_x)$ should be studied in the general case. For example, when does 
$W$ act on any of $H^*((\B_x)^S)$,  $H^*_S(\B_x)$ or even ${\iota^*}(H^*_S(\B))$?


\begin{thebibliography}{99}
\bibitem[A-H]{AH} Abe, H.; Horiguchi, T.
{\em The torus equivariant cohomology rings of Springer varieties}. Topology Appl. 208 (2016), 143--159. 

\bibitem[A-B]{AB} Atiyah, M.; Bott, R. {\em  The moment map and equivariant cohomology.} Topology 23 (1984), 1--28. 

\bibitem[Brion]{BRION} Brion, M. {\em Equivariant cohomology and equivariant intersection theory.} Notes by Alvaro Rittatore. NATO Adv. Sci. Inst. Ser. C Math. Phys. Sci., 514, Representation theories and algebraic geometry (Montreal, PX, 1997), 1--37, Kluwer Acad. Publ., Dordrecht, 1998.


\bibitem[JC1]{COMP} Carrell, J. B.  {\em Orbits of the Weyl group and a theorem of DeConcini and Procesi.} 
Compositio Math. 60 (1986), 45--52. 

\bibitem[JC2]{JBC} Carrell, J. B. {\em Torus actions and cohomology.} The adjoint representation and the adjoint action, 83--158, Encyclopaedia Math. Sci., 131, 
Springer, Berlin, 2002. 

\bibitem[D-L-P]{DLP} De Concini, C.; Lusztig, G.; Procesi, C. {\em Homology of the zero-set of a nilpotent vector field on a flag manifold.}
J. Amer. Math. Soc. 1 (1988), 15--34. 

\bibitem[D-P]{DP} De Concini, C.; Procesi, C. {\em Symmetric functions, conjugacy classes, and the flag variety}. Invent. Math., 64 (1981), 203--219.


\bibitem[G-McP]{GM} Goresky, M.; MacPherson, R. {\em On the spectrum of the equivariant cohomology ring}. Canad. J. Math. 62 (2010), 262--283.

\bibitem[H-S]{HS} Hotta, R.; Springer, T.A.  {\em A specialization theorem for certain Weyl group representations and 
an application to the Green polynomials of unitary groups}, Invent. Math. 41 (1977) 113--127.

\bibitem[Kraft]{KRAFT} Kraft, H. P. {\em Conjugacy classes and Weyl group representations}, Young tableaux and Schur functors in algebra and geometry (Torun,1980), 191�205, Ast�risque, 87�88, Soc. Math. France, Paris, 1981. 


\bibitem[K-P]{KP} Kumar, S.; Procesi, C.  {\em An algebro-geometric realization of equivariant cohomology
of some Springer fibers.}  Journal of Algebra 368 (2012), 70--74.


\bibitem[Ross]{Ross} Rossmann, W.  {\em Picard-Lefschetz theory for the coadjoint quotient of a semisimple Lie algebra.} (English summary) 
Invent. Math. 121 (1995), 531--578. 

\bibitem[Slo]{Slo} Slodowy, P.  {\em  Four lectures on simple groups and singularities.} Communications of the Math.
Inst., Rijksuniversiteit Utrecht, v. 11, ( 1 ) 1980.

\bibitem[Spa]{SPALT} Spaltenstein, N.  {\em The fixed point set of a unipotent transformation on the flag manifold.} 
Indag. Math. 38 (1976), 452--456. 

\bibitem[Spr1]{SPRINGER1} Springer, T. A. {\em Trigonometric sums, Green functions of finite groups and representations of Weyl groups.}  Invent. Math. 36 (1976), 173--207.

\bibitem[Spr2]{SPRINGER2}  Springer, T.  A. {\em A construction of representations of Weyl groups.} Invent. Math. 44 (1978),  279--293.

\bibitem[Treu]{Treu} Treumann, D.  {\em A topological approach to induction theorems in Springer theory.} 
Represent. Theory 13 (2009), 8--18. 


\end{thebibliography}
\end{document}